\newcommand{\be}{\begin{eqnarray}}
\newcommand{\ee}{\end{eqnarray}}
\newcommand{\nn}{\nonumber \\}

\newcommand{\CL}{{C\!\ell}}
\newcommand{\p}{\partial}
\newcommand{\B}[2]{B_{#1#2}}
\newcommand{\Bnb}[2]{B_{#1{\overline{#2}}}}
\newcommand{\Bbn}[2]{B_{{\overline{#1}}#2}}
\newcommand{\Bbb}[2]{B_{{\overline{#1}}\,{\overline{#2}}}}

\newcommand{\openone}{\mbox{1\kern -0.25em I}}
\newcommand{\openK}{\mbox{I\kern -0.25em K}}
\newcommand{\openZ}{\mbox{Z\kern -0.4em Z}}
\newcommand{\openR}{\mbox{I\kern -0.25em R}}
\newcommand{\openH}{\mbox{I\kern -0.25em H}}
\newcommand{\openM}{\mbox{I\kern -0.25em M}}
\newcommand{\openC}{\mbox{C\kern -0.55em I\hspace{0.25em}}}
\newcommand{\con}{\mbox{$\,$\rule{1ex}{0.4pt}\rule{0.4pt}{1ex}$\,$}}
\newcommand{\EOP}{\hfill\rule{4pt}{4pt}}
\newtheorem{dfn}{Definition}
\newtheorem{thrm}[dfn]{Theorem}
\documentstyle[a4]{article}

\begin{document}
\title{Hecke algebras as subalgebras of\\
       Clifford geometric algebras of multivectors}
\author{Bertfried Fauser\\
Eberhard-Karls-Universit\"{a}t\\
Institut f\"ur Theoretische Physik\\
Auf der Morgenstelle 14\\
72072 T\"ubingen {\bf Germany}\\
Electronic mail: Bertfried.Fauser@uni-tuebingen.de
}
\date{October 16, 1997}
\maketitle
\begin{abstract}
Clifford geometric algebras of multivectors are introduced which
exhibit a bilinear form which is not necessarily symmetric.
Looking at a subset of bi-vectors in $\CL(\openK^{2n},B)$, we
proof that theses elements generate the Hecke algebra
$H_{\openK}(n+1,q)$ if the bilinear form $B$ is chosen
appropriately. This shows, that $q$-quantization can be generated
by Clifford multivector objects which describe usually composite
entities. This contrasts current approaches which give deformed
versions of Clifford algebras by deforming the one-vector
variables. Our example shows, that it is not evident from a
mathematical point of view, that $q$-deformation is in any sense
more elementary than the undeformed structure.
\end{abstract}

\section{\protect\label{SEC-1}Introduction}

Recent developments in theoretical physics employ the so
called {\it noncommutative geometry} \cite{Con} or in a more
special case {\it $q$-deformed geometry}
\cite{Majid,Drinfeld,Jimbo}. The underlying structure is either
the $C$*-theory which incorporates also topological and
convergence aspects or else Hopf algebras, which model the
algebraic aspects of a theory \cite{Abe,Pittner}. It is
convenient to speak about $q$-symmetry since the spaces on which
$q$-symmetry acts tend to be {\it braided}. It is thus
convenient to study braided monoidal categories, e.g
\cite{Majid,Bautista,Ozi-multivectors} and many others. The main
idea is to introduce a braided tensor product algebra structure 
\be
(a \otimes b) (c \otimes d) &=& a \Psi(b\otimes c) d,
\ee
where $\Psi$ is a braiding. If $\Psi$ occurs to be trivial or
minus the flip operator $\Psi(a \otimes b) = -b\otimes a$, one
deals with the ordinary tensor product (bosons) or a
$\openZ_2$-graded version of it (fermions). A general braiding
leads thus to a general or braid statistics. The central
relations obeyed by braid groups are the Artin braid relations
\cite{Artin} 
\be
b_i b_{i+1} b_i &=& b_{i+1} b_i b_{i+1} \nn
b_i b_k &=& b_k b_i, \quad \vert i-k \vert \ge 2.
\ee
The first of them is actually equivalent to the so called
quantum Yang--Baxter equation, which is the special case of the
Yang--Baxter equation \cite{YBEQ} (in standard notation) 
\be
R_{12}(u) R_{13}(u+v) R_{23}(v) &=& R_{23}(v)R_{13}(u+v)R_{12}(v)
\ee
with the spectral parameters set to $v=u+v=u$.

There was a great progress in the theory of (quantum)
statistical mechanics, which originated in the development of
the inverse scattering method \cite{Faddeev} and the star
triangle relation \cite{ABF}, both methods having their roots
in braided symmetries, see e.g. \cite{YangGe}. There are lots of
models now solvable with this methods, Ising \cite{Ising} and
$N$-state Potts models \cite{Potts,FatZam}, Vertex
\cite{Baxter-buch} and IRF models \cite{Pasquier} may be
prominent examples. A further example might be given by the
(fractional) quantum Hall effect \cite{Bellisard}.
Furthermore, the unexpected connection between link invariants
and type III subfactors of von Neumann algebras unveiled by V.
Jones pushed low dimensional topology far ahead, see
\cite{Jones,Kauffmann-topology}. There is even a connection of
the Jones polynomial to quantum field theory \cite{Witten}.

A further branch of applications arises from the common believe,
that $q$-symmetry being more general than the usual bosonic or
fermionic ones and thus is more fundamental, see e.g.
\cite{Majid,Manin,WessZumino}. The natural thing to do
is thus to provide $q$-deformed versions of physical relevant
groups e.g. the Poincar\'e group \cite{Wess}. There is a strong
believe, that the fundamental constant $\hbar$ is involved in
this construction and that space-time should behave
``$q$-symmetric'' at small scales.  

The above mentioned situations when $q$-symmetry leads to
explicit results shear the feature of being effective or
composite models. There is no recent {\it evidence}\/ that
$q$-symmetry has to be used in fundamental interactions.
Moreover, it might be expected that a $q$-deformed Poincar\'e
group has an underlying structure which generates this symmetry. 

From a mathematical point, there is no harm in $q$-deforming all
structures which can be done so. But a physical application
requires an interpretation which seems currently not obvious,
but relays on rather abstract developments as quantum plains and
q-deformed or noncommutative geometry.

We are thus in the perplexing situation, that because of its
generality $q$-de\-for\-ma\-tion can be applied to nearly every
mathematical structure which is currently used in physics. But
we don't know in which cases it might be reasonable. To be able
to decide this question, it is a valuable advantage to have an
embedding of the mathematical structure which lays at the heart
of $q$-deformation, the Hecke algebra, in a larger framework.
From this outstanding point of view it might be possible to
decide if $q$-symmetry has to be applied to e.g. gravitation or
not. 

A very interesting approach to $q$-symmetry by spinors and
thereby also with the help of Clifford algebras can be found in
\cite{Bautista}. This approach, however, takes $q$-symmetry as
an elementary property. In the same spirit, the Clifford algebra
of a Hecke braid was constructed in \cite{Ozi-Hecke}.

We will provide a theorem, which shows that Hecke algebras, can
be obtained as subalgebras of certain Clifford algebras. This
subalgebras are generated by bi-vector elements and thus by
objects which are composed. Furthermore, since the
interpretation of Clifford algebraic expressions is well known,
we come to the end that $q$-symmetry is tightly interwoven with
{\it composite} structures, as was suggested already in
\cite{Fau-fkt}. This relation is seen from the fact that
$q$-symmetry is obtained in this approach as a multivector
symmetry. It is this relation, that makes us so suspicious
against a $q$-deformation of space-time as long as one does not
have a {\it microscopic description}\/ of these entities.
Hopefully our approach will open a possibility to clarify this
situation. 

\section{\protect\label{SEC-2}Clifford geometric algebra of
multivectors} 

There are many possibilities to introduce Clifford algebras,
each of them emphasize a different point of view. In our case,
it is of utmost importance to have the Clifford algebra build
over a graded linear space. This grading is obtained from the
space underlying a Grassmann algebra. The Clifford algebra is
then related to the endomorphism algebra of this Grassmann
algebra. This construction, the Chevalley deformation
\cite{Che}, was originally invented to be able to treat Clifford
algebras over fields of $char=2$, see appendix of \cite{MRiesz}
by Lounesto and \cite{Lou}. However, we use this construction in
an entire different context. With help of the construction of M.
Riesz \cite{MRiesz}, one is able to reconstruct the multivector
structure and thereby a correspondence between the linear spaces
underlying the Clifford algebra and the Grassmann algebra in use.
This reconstruction depends on an automorphism $J$, which is
arbitrary, see \cite{Fau-man}. In fact this is just the reversed
direction of our construction given below following Chevalley.

Let $T(V)$ be the tensor algebra build over the $\openK$-linear
space $V$. The field $\openK$ will be either $\openR$ or
$\openC$. With $V^0 \simeq \openK$ we have
\be
T(V) &=& \openK \oplus V \oplus V\otimes_{\openK} V \oplus 
\ldots . 
\ee
The tensor algebra is associative and unital. In $T(V)$ one has
bilateral or two-sided ideals, which can be used to construct
new algebras by factorization. As an example we define the
Grassmann algebra in this way.
\begin{dfn} The Grassmann algebra $\bigwedge(V)$ is the factor
algebra of the tensor algebra w.r.t. the bilateral ideal
\be
I_{Gr} &=& \{ y \mid y=a\otimes x \otimes x \otimes b,\quad
a,b\in T(V),\, x\in V\} \nn
\bigwedge(V) &=& \pi(T(V)) =\frac{T(V)}{I_{Gr}} =
\openK \oplus V \oplus V \wedge V \oplus \ldots .
\ee
The canonical projection $\pi : T(V) \mapsto \bigwedge(V)$ maps
the tensor product $\otimes$ onto the exterior or wedge product
denoted by $\wedge$.\EOP
\end{dfn}
One may note, that the factorization preserves the grading
naturally inherited by the tensor algebra, since the ideal
$I_{Gr}$ is homogeneous. Defining homogeneous parts of $T(V)$
and $\bigwedge(V)$ by $T^k(V)=V\otimes \ldots \otimes V$
and $\bigwedge^k(V)=V\wedge \ldots \wedge V$, $k$-factors, we 
obtain $\pi(T^k(V))=\bigwedge^k(V)$.

Proceeding to Clifford algebras requires a further structure,
the quadratic form.
\begin{dfn}\label{dfn-qform}
The map $Q : V \mapsto \openK$, satisfying ($\alpha\in \openK,
x,y\in V$)
\be\label{eq-qform}
i) && Q(\alpha x) = \alpha^2 Q(x) \nn
ii) && B_p(x,y) = \frac{1}{2}(Q(x+y)-Q(x)-Q(y)),
\ee
where $B_p(x,y)$ is a symmetric bilinear form is called a
quadratic form.\EOP
\end{dfn}
It is tempting to introduce an ideal $I_{\CL}$
\be\label{idealCL}
I_\CL&=&\{ y \mid y=a\otimes(x\otimes x -Q(x)\openone) \otimes
b, \quad a,b\in T(V),\, x\in V\}
\ee
to obtain the Clifford algebra by a factorization procedure.
However, since we are interested in arbitrary bilinear forms
underlying a Clifford algebra, we will take another approach,
which is much more reasonable for such a structure. Furthermore,
the Clifford algebra {\it does not}\/ have an intrinsic
multivector structure, but is {\it only}\/ $\openZ_2$ graded,
since the ideal $I_{\CL}$ is inhomogeneous. 

Let $V^*$ be the space of linear forms on $V$, i.e. $V^* \simeq
lin[V,\openK]$. Elements $\omega\in V^*$ act on elements $x\in
V$, but there is {\it no natural}\/ identification between $V$ and
$V^*$. However, we can find a set of $x_i$ which span $V$ and
dual elements $\omega_k$ acting on the $x_i$ in a {\it
canonical}\/ way
\be
\omega_k(x_i) &=& \delta_{ki}.
\ee
This allows to introduce a map $* : V \mapsto V^*$, $ x_i^* =
\omega_i$ which may be called Euclidean dual isomorphism
\cite{HaftSaller}. The two spaces $(V^*,V)$ connected by this
duality constitute a pairing $<.\mid .>\, : V^*\times V \mapsto
\openK$. $V^*$ are isomorphic to $V$ in finite dimensions, so it is natural to
build a Grassmann algebra $\bigwedge(V^*)$ over it. This is the
algebra of Grassmann multiforms.

It is further a natural thing to extend the pairing of the
grade-one space and its dual to the whole algebras
$\bigwedge(V)$ and $\bigwedge(V^*)$, as can be seen by its
frequent occurrence in literature
\cite{Fau-man,Ozi-Grass-Cliff,Fau-pos,Fau-diss,AblLou,Lou} and
others. This can be done by the 
\begin{dfn}\label{multiaction}
Let $\tau,\eta \in \bigwedge(V^*)$, $\omega \in V^*$, $u,v \in
\bigwedge(V)$ and $x\in V$, then we can define a canonical
action of $\bigwedge(V^*)$ on $\bigwedge(V)$ requiring
\be\label{tri-rel}
i) && \omega(x) = <\omega \mid x > \nn
ii) && \omega(u\wedge v) = w(u)\wedge v + \hat{u}\wedge
\omega(v) \nn
iii) && (\tau\wedge\eta)(u) = \tau(\eta(u))
\ee
where $\hat{u}$ is the main involution $\hat{V}=-V$ extended to
$\bigwedge(V)$. 
\EOP 
\end{dfn}
In fact we have given by definition \ref{multiaction} an
isomorphism between the Grassmann algebra of multiforms
$\bigwedge(V^*)$ and the dual Grassmann algebra
$[\bigwedge(V)]^*$. This can be made much clearer in writing 
\be&
y\con x = \omega_y(x) = < \omega_y \mid x > = B(y,x),
&\ee
where we have used the canonical identification of $V$ and
$V^*$ via the pairing. One should be very careful in the
distinction of $\bigwedge(V^*)$ and $[\bigwedge(V)]^*$, since
they are isomorphic but not equivalent. Furthermore, we
emphasize that in writing $y\con$ we make explicitly use of a
{\it special}\/ dual isomorphism encoded in the contraction
\be
&.\con : V \mapsto V^*& \nn
& y \rightarrow y\con = \omega_y.&
\ee
Since there is no natural, that is mathematically motivated or
even better functorial relation between $V$ and $V^*$, we are
called to seek for {\it physically motivated reasons}\/ to
select a pairing. This freedom will enable us in section
\ref{SEC-3} to give a proof of our main theorem.

\begin{thrm}
Let $(V,Q)$ be a pair of a $\openK$-linear space $V$ and $Q$ a
quadratic form $Q$ as defined in \ref{dfn-qform}. There exists
an injection $\gamma$ called Clifford map from $V$ into the
associative unital algebra $\CL(V,Q)$ which satisfies
\be\label{eq12}
\gamma_x \gamma_x &=& Q(x)\openone.
\ee
\EOP
\end{thrm}
\begin{dfn}\label{Cliffalg}
The (smallest) algebra $\CL(V,Q)$ generated by $\openone$ and
$\gamma_{x_i}$, $\{x_i\}$ span $V$ is called (the) Clifford
algebra of $Q$ over $V$. 
\EOP
\end{dfn}
By polarization of the relation (\ref{eq12}) we get the usual
commutation relations; $x,y\in V$
\be
\gamma_x \gamma_y + \gamma_y\gamma_x &=& 2B_p(x,y)\openone,
\ee
where $B_p(x,y)$ is the symmetric polar form of $Q$ as defined
in (\ref{eq-qform}).

\noindent {\bf Remarks:} i) We could have obtained this result
directly by factorization of the tensor algebra with the ideal
(\ref{idealCL}). ii) There exists Clifford algebras with are
universal, in this case it is convenient to speak from {\it
the}\/ Clifford algebra over $(V,Q)$. iii) If $V\simeq \openK^n
\simeq \openC^n$ or $\openR^n$, we denote $CL(V,Q)$ also by
$CL(\openC^n)\simeq\CL_n$ or $\CL(\openR_{p,q})$ where the pair
$p,q$ enumerate the number of positive and negative eigenvalues
of $Q$. We can as well give the dimension $n$ and signature
$s=p-q$ to classify all quadratic forms over $\openR$. In the
case of the complex field, on remains with the dimension as can
be seen e.g. from the Weyl unitary trick, letting $x_i
\rightarrow i x_i$ which flips the sign. We do not use
sesquilinear forms here, which could be included
nevertheless.\EOP 

We will now use Chevalley deformations to construct the Clifford
algebra of multivectors. The main idea is, that one can decompose
the Clifford map as
\be\label{decompose}
\gamma_x &=& x\con + x\wedge \, .
\ee
There is thus a natural action of $\gamma_x$ on $\bigwedge(V)$. 
\begin{thrm}[Chevalley] 
Let $\bigwedge(V)$ be the Grassmann algebra over $V$ and
$\gamma : V \mapsto End(\bigwedge(V))$, be defined as in
(\ref{decompose}), then $\gamma$ is a Clifford map.
\EOP 
\end{thrm}
We have shown that $\CL$ is a subalgebra of the endomorphism
algebra of $\bigwedge(V)$,
\be
\CL &\subseteq& \bigwedge(V).
\ee
It is possible to interpret $x\con$ and $x\wedge$ as
annihilating and creation {\it operators}\/ (on the Grassmann
algebra) \cite{Crum}.

With help of the relations (\ref{tri-rel}) we can then lift this
Clifford map to multivector actions. {\it No symmetry
requirement has to be made on the contraction.}\/ This leads to
the 
\begin{dfn}[Clifford algebra of multivectors]\label{def-Cliff}
Let $B : V\times V \mapsto \openK$ be an arbitrary bilinear form.
The Clifford algebra $\CL(V,B)$ obtained from lifting the
Clifford map 
\be&
\gamma_x = x\con + x\wedge = <x \mid .> + x\wedge = B(x,.)+x\wedge
&\ee
to $End(\bigwedge(V))$ using the relations (\ref{tri-rel}) is
called Clifford algebra of multivectors. \EOP
\end{dfn}
Note, that $B(x,.)=\omega_x$ is a map from $V \rightarrow V^*$
and {\it incorporates}\/ a dual isomorphism. It is clear from
the construction that $\CL(V,B)$ has a multivector structure or
say a $\openZ_n$-grading inherited from the Grassmann algebra
$\bigwedge(V)$.  

$B$ admits a decomposition into symmetric and antisymmetric
parts $B=G+F$. The symmetric part $G=B_p$ corresponds to a
quadratic form $Q$, see definition \ref{dfn-qform}.
\begin{thrm}
The Clifford algebra $\CL(V,Q)\simeq\CL(V,G)$ is isomorphic as
Clifford algebra to $\CL(V,B)$, if $B$ admits a decomposition
$B=G+F$, $G^T=G$, $F^T=-F$.\EOP
\end{thrm}
A proof can be found for low dimensions in \cite{AblLou} and in
general in \cite{Ozi-multivectors}. However, this result was
implicitly known to physicists, see \cite{Fau-ver,StuBor}. In
fact, this is the old Wick rule of QFT. We will insist on the
$\openZ_n$-grading and therefore carefully distinguish Clifford
algebras of multivectors with a common quadratic form $Q$ but
different contractions $B$. Only this generalization makes it
possible to find Hecke algebras as subalgebras in Clifford
algebras. 

We give some further notations. Let $\{j_i\}$ be a set of
elements spanning $V \simeq \,<j_1,\ldots,j_n>$ and $\{\p_k\}$
be a set of dual elements. Building the Grassmann algebras
$\bigwedge(V)$ and $\bigwedge(V^*)$ and defining the action of
the forms via (\ref{tri-rel}), one obtains the relations
\be
i) && j_i\wedge j_i = 0 = \p_i \wedge \p_i \nn
ii)&& \p_i j_k + j_k \p_i = B_{ik}+B_{ki} = 2G_{ik}.
\ee
The space ${\bf V}=V\oplus V^*$ is thus spanned by (note the
reversed order of indices for the $\p$ elements)
\be
\{e_1,\ldots, e_{2n} \} &=& \{j_1,\ldots,j_n,\p_n,\ldots,\p_1 \}.
\ee
To have a simple notation, we introduce bared indices $i\in
\{1,\ldots,n\}$ 
\be
e_{\overline{i}} &=& e_{2n+1-i},
\ee
or equivalently
\be
e_i=j_i &\quad& e_{\overline{i}}=e_{2n+1-i}=\p_i.
\ee
The contraction on $\CL({\bf V},B)$ is then written as
\be
[B(e_i,e_j)]=[B_{ij}]&=&
\begin{array}{|cc|}
B^{jj}_{uv} & B^{j\p}_{st} \\
B^{\p j}_{rq} & B^{\p\p}_{xy}
\end{array} \Leftrightarrow
\begin{array}{|cc|}
[M_{rs}]^{jj} & [\Bnb tu^1]^{j\p} \\
{}[\Bbn vw^2]^{\p j} & [N_{\overline{x}\,\overline{y}}]^{\p\p}
\end{array},
\ee
where the overscripts indicate the type of the base element.
Indices of blocks run in $\{1\ldots n\}$. Note, that the
matrices $B^1,B^2$ and $N$ are not directly submatrices of $B$
because of our bared index notation. Introducing a $n\times n$
matrix $J=\delta_{i,n+1-k}$ we can identify them as
\be
& B^{j\p} = B^1J,\quad B^{\p j} = JB^2,\quad B^{\p\p}=JNJ.&
\ee
We could handle the $2n$ dimensional complex case as
$\CL(\openR^{2n+1},B)$, but we will restrict ourself to the even
dimensional case and look at $\CL(\openC^{2n}) \simeq
\openC\otimes \CL(\openR^{2n})$ as a complexification.

\section{\protect\label{SEC-3}Hecke algebra as bi-vector subalgebra} 

\begin{dfn}\label{dfn-Hecke}
The Hecke algebra $H_{\openK}(n+1,q)$ has the following
presentation  
\be\label{def-Hecke}
&
\begin{array}{lll}
i)    & b_i^2=(1-q) b_i +q\openone  & \mbox{Hecke condition} \\
ii)   & b_i b_k = b_k b_i \quad \vert i-k \vert \ge 2 & 
\mbox{commutator}\\
iii)  & b_i b_{i+1} b_i = b_{i+1} b_i b_{i+1} &
\mbox{Artin braid relation \cite{Artin}}
\end{array}&
\ee
with generators $\openone,b_i$, $i\in \{ 1,\ldots,n\}$, see
\cite{Bour}. \EOP
\end{dfn}
Our goal is to find an identification of the $b_i$ generators
as bi-vectors in an appropriate Clifford algebra
$\CL(\openR^{2n},B)$ or $\CL(\openC^{2n},B)$. We can formulate
our results in the following\\ 
\begin{thrm}\label{mainthrm}
The Hecke algebra $H_{\openK}(n+1,q)$ of definiton
(\ref{dfn-Hecke}) is a subalgebra of the Clifford algebra
$\CL(\openK^{2n},B)$ of definition \ref{def-Cliff} with the 
following identifications:
\be\label{def-bi}
i) && b_i:=e_i \wedge e_{\overline{i}}= 
e_i \wedge e_{2n+1-i} \equiv j_i \wedge \partial_i
\quad i\in \{1,\ldots,n\} \nn
ii) && B:= [\B ij ] = \begin{array}{|cc|}
[\B rs]^{jj} & [\B tu ]^{j\p} \\
{}[\B vw ]^{\p j} & [\B xy]^{\p\p}
\end{array}
\ee
where the submatrices of $B$ satisfy the conditions
\be\label{formofb}
B^{jj} \equiv M_{rs} &=& \frac{1}{2}(M_{rs}-M_{sr}) \nn
JB^{\p\p}J \equiv N_{\overline{x}\,\overline{y}} &=&
\frac{1}{2}(N_{\overline{x}\,\overline{y}}
          -N_{\overline{y}\,\overline{x}})\nn 
B^{j\p}J \equiv [\Bnb tu^1 ] &=&
[\Bnb tu +(q-\Bnb tu)\delta_{t,\overline{u}}] \nn
JB^{\p j} \equiv [\Bbn vw^2 ] &=&
[-\Bnb wv +(1+q)\delta_{w,\overline{v}} 
  +q\delta_{w+1,\overline{v}}
  + \delta_{w,\overline{v}+1} ]
\ee
\EOP
\end{thrm}

Proof: We determine the constraints on the bilinear form $B$ by
a direct calculation of the consequences of the relations
(\ref{def-Hecke}). 

\noindent{\bf i)} We try to identify the bi-vector elements
$b_i$ from (\ref{def-bi}) with Hecke generators also denoted by
$b_i$ from (\ref{def-Hecke}). Since we insist on the multivector
structure inherited from the Grassmann multivectors underlying
the Clifford multivectors, we have to fulfill in any case the
condition 
\be\label{NB1}
e_i e_i &=& e_i \wedge e_i = B_{ii} =0, \quad (\Bbb ii =0).
\ee
The Hecke relation (\ref{def-Hecke}-i) leads with
\be
b_i &=& j_i \wedge \p_i = j_i \p_i -\Bnb ii 
\ee
to
\be\label{calc-i}
b_i^2 &=& (j_i \wedge \p_i )^2 = (j_i \p_i -\Bnb ii ) j_i \wedge
\p_i \nn
   &=& j_i [\Bbn ii \p_i - j_i \p_i^2] -\Bnb ii j_i \wedge \p_i
\nn
   &=& \Bbn ii \Bnb ii -(\Bbn ii -\Bnb ii ) j_i \wedge \p_i \nn
   &=& \Bbn ii \Bnb ii -(\Bbn ii -\Bnb ii ) b_i \nn
   &\stackrel{\displaystyle !}{\displaystyle =}&
   (1-q) b_i +q.
\ee
We get as solutions
\be\label{NB2}
\Bnb ii  &=& q \quad \mbox{or} \quad -1 \nn
\Bbn ii &=& 1 \quad \mbox{or} \quad -q.
\ee
We will chose $\Bnb ii =q$, $\Bbn ii =1$. The overall minus sign
does not matter in our considerations. Including the nilpotency
of the Grassmann sources $j$ and $\p$, we obtained $4n$
constrains on $B$.

\noindent{\bf ii)} The commutator relation (\ref{def-Hecke}-ii),
which is valid for $\vert k-i\vert \ge 2$, can be calculated
along the same lines as in (\ref{calc-i}). This results in
\be
b_i b_k - b_k b_i &=&
\big(
\Bbn ik \Bnb ik -\Bbn ki \Bnb ki -\Bbb ik \B ik +\Bbb ki \B ki
\big)\openone \nn
&&+\big(\Bbn ik + \Bnb ki) j_i\wedge \p_k 
      -(\Bnb ik + \Bnb ki\big) j_k\wedge \p_i \nn
&&-\big(\Bbb ik + \Bbb ki)j_i \wedge j_k 
      -(\B ik + \B ki\big) \p_i \wedge \p_k \nn
&\stackrel{\displaystyle !}{\displaystyle =}& 0.
\ee
Therefrom we obtain
\be\label{NB3}
\B ik &=& - \B ki \nn
\Bbb ik &=& -\Bbb ki \nn
\Bnb ik &=& - \Bbn ki \nn
\ee
if $\vert i-k \vert \ge 2$. This leads to 3n(n-2)/2 constrains on
$B$. 

\noindent{\bf iii)} The third relation is somewhat lengthy and
yields
\be
b_i b_{i+1} b_i - b_{i+1} b_i b_{i+1} &=&  
   (1+q)\big(
        \B i{i+1}\Bbb {i+1}i -\B {i+1}i \Bbb i{i+1}\big)
        \openone \nn &&
+\big(
      (\Bnb {i+1}i+\Bbn i{i+1})(\Bnb i{i+1}+\Bbn {i+1}i)
      \nn && ~~
     -(\B i{i+1}+\B {i+1}i)(\Bbb {i+1}i +\Bbb i{i+1}) 
  -q\big) 
j_i \wedge \p_{i} \nn &&
+\big(
      (\Bbb i{i+1}+\Bbb {i+1}i)(\B {i+1}i +\B i{i+1})
      \nn && ~~
     -(\Bbn {i+1}i+\Bnb i{i+1})(\Bnb {i+1}i+\Bbn i{i+1})
  +q\big)
j_{i+1} \wedge \p_{i+1} \nn &&
- (1+q)(\B i{i+1}+ \B {i+1}i) \p_{i} \wedge \p_{i+1} \nn &&
+ (1+q)(\Bbb i{i+1}+ \Bbb {i+1}i) j_i \wedge j_{i+1} \nn 
&\stackrel{\displaystyle !}{\displaystyle =}& 0.
\ee
This leads to
\be\label{NB4}
\B {i+1}i &=& -\B i{i+1} \nn
\Bbb i{i+1} &=& -\Bbb {i+1}i \nn
\Bbn i{i+1} &=& 1- \Bnb {i+1}i \nn
\Bbn {i+1}i &=& q- \Bnb i{i+1},
\ee
which are $4(n-1)$ further constrains on $B$. All in all, we
have to impose the constrains given in
(\ref{NB1},\ref{NB2},\ref{NB3}) and (\ref{NB4}) on the bilinear
form $B$ of $4n^2$ arbitrary $\openK$-valued parameters. We
obtain 
\be
\# \mbox{constrains} &=& \frac{3n^2+13n-8}{2}
\ee
and remain with 
\be
\# \mbox{degrees of freedom} &=&
\frac{5n^2-13n+8}{2}.
\ee
The explicit form of $B$ can be derived from the constraints to
be of the form (\ref{formofb}). Since we remain with superfluous
degrees of freedom, which might be set to zero, we have derived a
whole set of Hecke algebra embeddings in $\CL(V,B)$.\EOP

Since we have proofed theorem \ref{mainthrm} for the quadratic
Hecke condition (\ref{def-Hecke}-i), we have to give some
comments on other choices of the quadratic or higher relations.
One finds in literature at least the following types of relations 
\be
b_i^2 &=& (1-q)b_i +q \nn
t_i^2 &=& (q-q^{-1}) t_i +1 \nn
e_i^2 &=& \tau e_i \nn
u_i^{Q} &=& 1 .
\ee
In general one has a quadratic --or higher order, see $u_i$--
relation 
\be
x_i^2 &=& g(q) x_i + h(q)
\ee
where $g$ and $h$ are mereomorphic functions of $q$. The
question, if there is a transformation connecting the $b_i$'s
and in general the $x_i$'s addresses the number of equivalence
classes of quadratic relations. We can reformulate the above
equations into $Q(x)=E$ where $E$ is a constant, and have to
classify quadratic forms. This can be done over $\openR$ and
$\openC$ with help of the Brauer group $B(\openK)$,
\cite{Hahn,Lou}. This group is trivial for $\openC$ since the
complex field is algebraically closed and isomorphic to
$\{-1,+1\}$ as a multiplicative group in the case of $\openR$
\be
B(\openK) &\cong& \frac{\openK^*}{\openK^2} .
\ee
However, this classification does only take (\ref{def-Hecke}-i)
into account. It is easy to calculate, that a transformation of
the type
\be
x_i &=& a(q) + b(q) b_i
\ee
with mereomorphic $a(q),b(q)$ does not alter
(\ref{def-Hecke}-ii). But (\ref{def-Hecke}-iii) is in general
{\it not} invariant under such a transformation, which is well
known in literature. As an example, one arrives at the relations
of a Temperly-Lieb algebra \cite{TemLie}, where the third
relation is given as, $e_i = (q\openone+b_i)/(q+1)$,
$\tau=(2+q+q^{-1})^{-1}$ compare (\ref{tri-rel}-iii)
\be
e_i e_{i+1} e_i -\tau e_i &=&
e_{i+1} e_i e_{i+1} -\tau e_{i+1}.
\ee
Such a relation can easily be obtained in our approach, simply
by another choice of the bilinaear form $B$, or of course by the
above transformation. It remains however to find a connection
between traces employed in the phenomenological models and
states on our algebra. Such stats were introduced in
\cite{Fau-vac} and provide a very rigid structure on $\CL(V,B)$.
These states are necessary to be able to calculate invariants
and physical outcomes and to be able to show (in)equivalence
between different presentations. This intriguing problem will be
addressed elsewhere. 

\section{\protect\label{SEC-4}Conclusion}

Our main tool were Clifford geometric algebras of multivectors,
which provide a generalization of ordinary Clifford algebras of
quadratic spaces to such of a pair ($V,B$). The bilinear form
does not having any symmetry i.g. and is thus not bijectively
connected to a quadratic form. Every bilinear form $B$ with the
same symmetric part $G$ gives rise to the same Clifford algebra.
Taking the full $B$ into account allows one to endow the
Clifford algebra with a unique $\openZ_n$-grading. The Clifford
algebra corresponding to $B$ build over the $\openZ_n$-graded
space $\bigwedge(V)$ is called Clifford algebra of multivectors
\cite{Ozi-multivectors}. 

We proofed the theorem, that due to an appropriate choice of the
bilinear form $B$ in $\CL(\openK^{2n},B)$, it is possible to
find $n$ bivectors $b_i$ which generate the Hecke algebra
$H_{\openK}(n+1,q)$. The proof was straight forward. Since we
got a large number of remaining freedoms in the Clifford
bilinearform $B$, this parameters might be used to have spectral
parameters in the braid relation, which then mutates into the
Yang-Baxter equation. This will be considered elsewhere.

Since the Clifford algebra has already an interpretation in
physical terms, we have to look at the $b_i$ generators as
composite entities. This supports our opinion stated in the
introduction and also promoted in \cite{Fau-fkt} that
$q$-symmetry might be connected with composite effects. A
decision of this conjecture requires further work, especially on
the states involved in the calculation of invariants, see
\cite{Fau-vac}. 

\section{Acknowledgements}

This paper was inspired by the possibility to play with
``Clifford'' Ver. 3, a Maple V Release 4 add-on programmed by
Prof. Rafal Ablamowicz, Gannon Erie \cite{Abla}. His generous
help in extending the features, collected in the Cli3plus bonus
pack, were a valuable help.
The support of Mrs. Ursula Wieland is gratefully acknowledged. 

\end{document}